\documentclass[12pt,a4paper]{article}
\usepackage{amsmath,amsfonts}
\usepackage[all]{xy}
\begin{document}

\begin{center}
{\bf Forms of Affine Space}\\
by\\
S. Subramanian\\[5mm]
\end{center}

Let $k$ be a field and $k^s$ its separable closure.  Let $A$ be a finitely
generated $k$ algebra such that $A\otimes_k k^s$ is isomorphic to the 
polynomial algebra in $n$-variables over $k^s$.  Let $V$ be an $n$ 
dimensional
 $k^s$ vector space so that 
$$
\varphi : \bigoplus_{n \geq 0} S^n V \rightarrow A\otimes_k k^s
$$
is an isomorphism of $k^s$ algebras.\\
(Here we adopt the convention that $S^0 V$ is the field $k^s$).

Let $V_0$ be a $k$-vector space of dimension $n$ contained in $V$ such that
$V_0 \otimes_k k^s$ is isomorphic to $V$ as a vector space via the 
inclusion 
map $V_0 \subset V$.  Then the $k$-algebra
$$
\bigoplus_{n\geq 0} S^n V_0
$$
(again with the conventions that $S^0 V =k$), is a $k$-subalgebra of $
{\displaystyle\bigoplus_{n\geq 0}} S^n V$. The map $\varphi$ induces a map of 
$k$-algebras,
$$
\overline{\varphi}: \bigoplus_{n\geq 0} S^n V_0 \rightarrow A \otimes_k k^s.
$$
It follows that there is a map
$$\psi : \bigoplus_{n\geq 0} S^n V_0 \rightarrow A\otimes_k 1 \subset A 
\otimes_k k^s
$$
of $k$-algebras and a factoring

$$
\xymatrix{
\displaystyle{\bigoplus_{n \geq 0}} S^nV_0  \ar[rd]_{\psi}
&\ar[r]^{\overline{\varphi}}&& A \otimes_k k^s\\
&&A \otimes 1 \ar[ru] &
}
$$

(We note that this factoring exists precisely because $k^s/k$ is a separable 
extension. In the inseparable case, this is not true).

>From the construction, it is obvious that after tensoring with $k^s, \ \psi 
\otimes k^s$ is isomorphic to $\varphi$.  Since $\varphi$ is an isomorphism of
 $k^s$ algebras, it follows that $\psi$ is an isomorphism of $k$-algebras.  We
 have thus shown:

\paragraph*{Theorem:} Let $k$ be any field and $k^s$ its separable 
closure.  
Let $X$ be an affine variety over $k$ which is isomorphic to affine $n$-space 
over the field extension $k^s$.  Then $X$ is isomorphic to affine $n$ space 
over $k$.

\vspace*{10mm}

{\bf References}
\begin{enumerate}
\item P. Russell, Some Formal Aspect of the Theorems of Manifold-Ramanujan, 
Proceedings of the International Colloquium on Algebra, Arithmetic and 
Geometry 2000, p 557-584, TIFR Mumbai, Narosa Publishing House.
\item T. Kambayashi, On the absence of nontrivial separable forms of the 
affine plane, Journal of Algebra, Volume  35, 1975, p. 449-456.
\item I.R. Safarerich, On Some Infinite Dimensional Groups, p. 208-212, Atti. 
Simposio Internaz. di Geometria Algebrica Roma, 1965.
\end{enumerate}
\end{document}